\documentclass[12pt]{article}
\usepackage{amssymb}
\usepackage{amsmath,amsthm}
\usepackage[latin1]{inputenc}
\usepackage{citesort}
\usepackage{graphicx}
 \newtheorem{remark}{Remark}

 \newtheorem{theorem}[remark]{Theorem}

\title{On the global offensive alliance number of a graph}

\author{J. M. Sigarreta\footnote{e-mail:\mbox{\tt
    josemaria.sigarreta\@@uc3m.es}}\\
{\em Department of Mathematics }\\ University Carlos III of Madrid\\
Avda. de la Universidad 30, 28911
Leganés (Madrid),  Spain \\
J. A. Rodr\'{\i}guez\footnote{e-mail:\mbox{\tt
    juanalberto.rodriguez\@@urv.net}} \\
{\em Department of Computer Engineering and Mathematics}\\
Rovira i Virgili University of Tarragona\\ Av. Pa\"{\i}sos Catalans
26, 43007 Tarragona, Spain}

\date{}

\begin{document}

\maketitle

\begin{abstract}
An offensive alliance  in a graph $\Gamma=(V,E)$ is a set of
vertices $S\subset V$ where for every vertex $v$ in its boundary it
holds that the majority of vertices in $v$'s  closed neighborhood
are in $S$. In the case of strong offensive alliance, strict
majority is required. An alliance $S$ is called global if it affects
every vertex in $V\backslash S$, that is, $S$ is a dominating set of
$\Gamma$. The offensive alliance number $a_o(\Gamma)$ (respectively,
strong offensive alliance number   $a_{\hat{o}}(\Gamma)$) is the
minimum cardinality of an offensive (respectively, strong offensive)
alliance in $\Gamma$. The global offensive alliance number
$\gamma_o(\Gamma)$   and the global strong offensive alliance number
$\gamma_{\hat{o}}(\Gamma)$ are defined similarly. Clearly,
$a_o(\Gamma)\le \gamma_o(\Gamma)$ and $a_{\hat{o}}(\Gamma)\le
\gamma_{\hat{o}}(\Gamma)$. It was shown in [{\it Discuss. Math.
Graph Theory} {\bf 24} (2004), no. 2, 263-275] that $ a_o(\Gamma)\le
\frac{2n}{3}$ and $ a_{\hat{o}}(\Gamma)\le \frac{5n}{6}$, where $n$
denotes the order of $\Gamma$. In this paper we obtain several tight
bounds on $\gamma_o(\Gamma)$ and $\gamma_{\hat{o}}(\Gamma)$ in terms
of several parameters of $\Gamma$. For instance, we show that
$\frac{2m+n}{3\Delta+1} \le \gamma_o(\Gamma)\le \frac{2n}{3}$ and
$\frac{2(m+n)}{3\Delta+2} \le\gamma_{\hat{o}}(\Gamma)\le
\frac{5n}{6}$, where $m$ denotes the size of $\Gamma$ and $\Delta$
its maximum degree (the last upper bound holds true for all $\Gamma$
with minimum degree greatest or equal to two).

\end{abstract}



\section{Introduction}

The study of defensive alliances in graphs, together with a variety
of other kinds of allian\-ces,  was introduced in
\cite{alliancesOne}. In the referred paper was initiated the study
of the mathe\-matical properties of alliances. In particular,
several bounds on the defensive alliance number were given. The
particular case of global (strong) defensive alliance was
investigated in \cite{GlobalalliancesOne}.

The study of offensive alliances was initiated by Favaron et al. in
\cite{favaron} where were derived several bounds on the offensive
alliance number and the strong offensive alliance number. On the
other hand, in  \cite{spectral} were obtained several tight bounds
on different types of allian\-ce numbers of a graph: (global)
defensive alliance number, global offensive alliance number and
global dual alliance number. In particular, was investigated the
relationship between the alliance numbers of a graph and its
algebraic connectivity, its spectral radius, and its Laplacian
spectral radius. A particular study of the alliance numbers, for the
case of planar graphs, can be found in \cite{planar}. Moreover,  for
the study of defensive alliances in the line graph of a simple graph
we cite \cite{linedefensive}.

The aim of this paper is to study mathematical properties of the
global  offensive alliance number and the global strong offensive
alliance number of a graph. We begin by stating some notation and
terminology. In this paper $\Gamma=(V,E)$ denotes a simple graph of
order $n$ and size $m$. The  degree of a vertex $v\in V$ will be
denoted by $\delta(v)$, the minimum degree will be denoted by
$\delta$, and the maximum degree by $\Delta$. The subgraph induced
by a set $S\subset V$ will be denoted by $\langle S\rangle$. For a
non-empty subset $S\subset V$, and a vertex $v\in V$, we denote by
$N_S(v)$ the set of neighbors $v$ has in $S$: $N_S(v):=\{u\in S:
u\sim v\}$.
 Similarly, we denote by
$N_{V\setminus S}(v)$ the set of neighbors $v$ has in $V\setminus
S$: $N_{V\setminus S}(v):=\{u\in V\setminus S: u\sim v\}$.  The
boundary of a set $S\subset V$ is defined as $\partial
(S):=\displaystyle\bigcup_{v\in S}N_{V\setminus S}(v).$

A non-empty set of vertices $S\subset V$ is called {\em offensive
alliance} if and only if for every $v\in \partial (S)$, $| N_S(v) |
\ge | N_{V\setminus S}(v)|+1.$ That is, a non-empty set of vertices
$S\subset V$ is called  offensive alliance if and only if for every
$v\in \partial (S)$, $2| N_S(v) | \ge \delta(v)+1.$

An offensive alliance $S$ is called {\em strong} if for every vertex
$v\in \partial (S)$, $| N_S(v) | \ge | N_{V\setminus S}(v)|+2.$ In
other words, an offensive alliance $S$ is called  strong if for
every vertex $v\in \partial (S)$, $2| N_S(v) | \ge \delta(v)+2.$

The {\em   offensive alliance number} (respectively, \emph{strong
offensive alliance number}), denoted  $a_{o}(\Gamma)$ (respectively,
 $a_{\hat{o}}(\Gamma)$), is defined as
the minimum cardinality of an  offensive alliance (respectively,
strong offensive alliance) in $\Gamma$.

 A non-empty set of vertices $S\subset V$ is a {\em global offensive
alliance} if for every vertex $v\in V\setminus S$, $| N_S(v) | \ge |
N_{V\setminus S}(v)|+1$. Thus, global offensive alliances are also
dominating sets, and one can define the {\em global offensive
alliance number}, denoted $\gamma_{o}(\Gamma)$, to equal the minimum
cardinality of a global offensive alliance in $\Gamma$. Analogously,
$S\subset V$ is a {\em global strong offensive alliance} if for
every vertex $v\in V\setminus S$, $| N_S(v) | \ge | N_{V\setminus
S}(v)|+2,$ and the  {\em global strong offensive alliance number},
denoted $\gamma_{\hat{o}}(\Gamma)$, is defined as the minimum
cardinality of a global strong offensive alliance in $\Gamma$.

In this paper we obtain several tight bounds on $\gamma_o(\Gamma)$
and $\gamma_{\hat{o}}(\Gamma)$ in terms of several parameters of
$\Gamma$. For instance, we show that
\begin{equation}
\left\lceil\frac{2m+n}{3\Delta+1}\right\rceil \le
\gamma_o(\Gamma)\le \left\lfloor\frac{2n}{3}\right\rfloor
\end{equation}
and
\begin{equation}
\left\lceil\frac{2(m+n)}{3\Delta+2}\right\rceil
\le\gamma_{\hat{o}}(\Gamma)\le \left\lfloor\frac{5n}{6}\right\rfloor
\end{equation}
(the
last upper bound holds true for all $\Gamma$ with minimum degree
greatest or equal to two).

\section{Bounding above the global offensive alliance number}

It was shown in \cite{favaron} that the offensive alliance number of
a graph of order $n\ge 2$ is bounded by
\begin{equation} \label{fab1}
a_o(\Gamma)\le \left\lfloor\frac{2n}{3}\right\rfloor,  \quad
a_o(\Gamma) \le \left\lfloor\displaystyle\frac{\gamma(\Gamma) +
n}{2}\right\rfloor,
\end{equation}
where $\gamma(\Gamma)$ denotes de domination number of $\Gamma$,
 and the strong offensive alliance number of a graph of order $n\ge 3$ is bounded
by
\begin{equation} \label{fab2}
a_{\hat{o}}(\Gamma)\le \left\lfloor\frac{5n}{6}\right\rfloor.
\end{equation}
%It is clear that .
Clearly, $a_o(\Gamma)\le \gamma_o(\Gamma)$ and
$a_{\hat{o}}(\Gamma)\le \gamma_{\hat{o}}(\Gamma)$.  Now we are going
to obtain the above bounds for the case of global alliances.

\begin{theorem}
For all connected graph $\Gamma$ of order $n\ge 2$,
\begin{itemize}
\item[i{\rm )}] $\gamma_o(\Gamma)\le \min\left\{n-\alpha(\Gamma),
\left\lfloor\displaystyle\frac{n+\alpha(\Gamma)}{2}\right\rfloor\right\}$,
 where $\alpha(\Gamma)$ denotes the independence number of $\Gamma$;
\item[ii{\rm )}] $\gamma_o(\Gamma)\le \left\lfloor\displaystyle\frac{2n}{3}\right\rfloor$;
\item[iii{\rm )}] $\gamma_o(\Gamma) \le \left\lfloor\displaystyle\frac{\gamma(\Gamma) +
n}{2}\right\rfloor$, where $\gamma(\Gamma)$ denotes the domination
number of $\Gamma$;
\item[iv{\rm )}] $\gamma_o(\Gamma) \le  \left\lfloor\displaystyle\frac{n(2\mu-\delta)}{2\mu}\right\rfloor$,
where $\mu$ denotes the Laplacian spectral radius of $\Gamma$ and
$\delta$ denotes its minimum degree.
\end{itemize}

\end{theorem}

\begin{proof}
Let $S\subset V$ be an independent set of maximum cardinality
$\alpha(\Gamma)$. Since the set $V\backslash S$ is a global
offensive alliance in $\Gamma=(V,E)$, then
\begin{equation} \label{eq1}
\gamma_o(\Gamma)+\alpha(\Gamma)\le n.
 \end{equation}
If $|V\backslash S|=1$, then $\Gamma=K_{1,n-1}$ and
$\gamma_o(\Gamma)=1$. If $|V\backslash S|\neq 1$, let $V\backslash
S=X\cup Y$ be a partition of $V\backslash S$ such that the edge-cut
between $X$ and $Y$ has the maximum cardinality. Suppose  $|X|\le
|Y|$. For every $v\in Y$, $|N_S(v)|\ge 1$ and $|N_X(v)|\ge
|N_Y(v)|$. Therefore, the set $W=S\cup X$ is a global offensive
alliance in $\Gamma$, i.e., for every $v\in Y$, $|N_W(v)|\ge
|N_Y(v)|+1$. Then we have, $2|X|+\alpha(\Gamma)\le n$ and
$\gamma_o(\Gamma)\le|X|+\alpha(\Gamma) $.
 Thus,
 \begin{equation}\label{eq2}
2\gamma_o(\Gamma) -\alpha(\Gamma)\le n.
\end{equation}
The bounds {\it i}) and {\it ii}) follow from (\ref{eq1}) and
(\ref{eq2}).

The proof of {\it iii}) follows in the spirit of the proof of
(\ref{eq2}): in this case we take $S\subset V$  as a dominating set
of minimum cardinality. Finally, it was shown in \cite{partition}
that $$\alpha(\Gamma)\le \frac{n(\mu-\delta)}{\mu}.$$ Thus, by
(\ref{eq2}) we obtain {\it iv}).
\end{proof}

The above bounds are attained, for instance, for the cocktail-party
graph $\Gamma=K_6-F\cong K_{2,2,2}$ where $n=\mu=6$, $\delta=4$,
$\alpha(\Gamma)=\gamma(\Gamma)=2$ and $\gamma_o(\Gamma)$=4.

In the spirit of the proof of {\it iii}) we obtain
\begin{equation}\label{eq6}
2\gamma_o(\Gamma)-\gamma_{c}\le n,
\end{equation}
 where $\gamma_c(\Gamma)$ denotes the connected-domination number of $\Gamma$.
 Moreover, it was shown in  \cite{cconex}  that if  $\Gamma$ is a connected graph of order  $n$
 and maximum degree $\Delta$, then
\begin{equation}\label{eq9}
\gamma_{c}\le n-\Delta.
\end{equation}
Thus, by  (\ref{eq6}) and (\ref{eq9}) we obtain
\begin{equation}
\gamma_o(\Gamma)\le \left\lfloor \frac{2n-\Delta}{2}\right\rfloor.
\end{equation}
This bound improves {\it ii}) if $\Delta > \frac{2n}{3}$.

\begin{theorem} \label{strongBelow}
For all connected graph $\Gamma$ of order $n$,
\begin{itemize}
\item[i{\rm )}] $\gamma_{\hat{o}}(\Gamma)\le
\left\lfloor\displaystyle\frac{n+\gamma_2(\Gamma)}{2}\right\rfloor$,
where $\gamma_2(\Gamma)$ denotes the $2$-domination number of
$\Gamma$.
\end{itemize}
 If the minimum degree of $\Gamma$ is greatest or equal to
two, then
\begin{itemize}
\item[ii{\rm )}]
$\gamma_{\hat{o}}(\Gamma)\le n-\alpha(\Gamma)$, where
$\alpha(\Gamma)$ denotes the independence number of $\Gamma$;
\item[iii{\rm )}] $\gamma_{\hat{o}}(\Gamma)\le
 \left\lfloor\displaystyle\frac{5n}{6}\right\rfloor$;
\item[iv{\rm )}] if $\Gamma$ is a cubic
graph, then $\gamma_{\hat{0}}(\Gamma)\leq \left\lfloor
\displaystyle\frac{3n}{4}\right\rfloor.$
\end{itemize}
\end{theorem}

\begin{proof}

Let $H\subset V$ be a 2-dominating set of minimum cardinality. If
$|V\backslash H|=1$, then $\gamma_2(\Gamma)=n-1$ and $\gamma_{\hat
o}(\Gamma)\le n-1$. If $|V\backslash H|\neq 1$, let $V\backslash
H=X\cup Y$ be a partition of $V\backslash H$ such that the edge-cut
between $X$ and $Y$ has the maximum cardinality. Suppose  $|X|\le
|Y|$. For every $v\in Y$, $|N_H(v)|\ge 2$ and $|N_X(v)|\ge
|N_Y(v)|$. Therefore, the set $W=H\cup X$ is a global strong
offensive alliance in $\Gamma$, i.e., for every $v\in Y$,
$|N_W(v)|\ge |N_Y(v)|+2$. Then we have,
\begin{equation}\label{eq3} 2|X|+\gamma_2(\Gamma)\le n\end{equation} and
\begin{equation}\label{eq4} \gamma_{\hat{o}}(\Gamma)\le|X|+\gamma_2(\Gamma) .\end{equation}
 Thus, by (\ref{eq3}) and (\ref{eq4}), {\it i}) follows.

Let $S\subset V$ be an independent set of maximum cardinality
$\alpha(\Gamma)$. Since $\delta \ge 2$, the set $V\backslash S$ is a
global strong  offensive alliance in $\Gamma=(V,E)$. Hence, {\it
ii}) follows. On the other hand, it was shown in \cite{Cockayne}
that
 \begin{equation}\label{eq5} \delta\ge
 2\Rightarrow\gamma_2(\Gamma)\le \frac{2n}{3}.
 \end{equation}
  So, by  {\it  i}) and (\ref{eq5}),  {\it iii}) follows.

  Finally, if $\Gamma$ is connected with maximum degree
$\Delta \leq 3$, then for all  global strong offensive alliance $S$
such that $|S|=\gamma_{\hat{0}}(\Gamma)$, $V\backslash S$ is an
independent set. Thus, $m\leq
3(n-\gamma_{\hat{0}}(\Gamma))+\gamma_{\hat{0}}(\Gamma)$. Hence, the
result follows.
\end{proof}

%XXXXXXXXXXXXXXXXX ofensiva XXXXXXXXXXXXXXXXXXXx
%\begin{figure}[h]
%\begin{center}
%\caption{ } \label{fig1}
%\includegraphics[width=0.3\textwidth]{ofensivateoiii}
%\includegraphics[width=0.3\textwidth]{ofensivateo3ab}
%\end{center}
%\end{figure}

\begin{figure}[h]
\begin{center}
\caption{} \label{fig1} %\vspace{-1,5cm}
\includegraphics[angle=0, width=3.50cm]{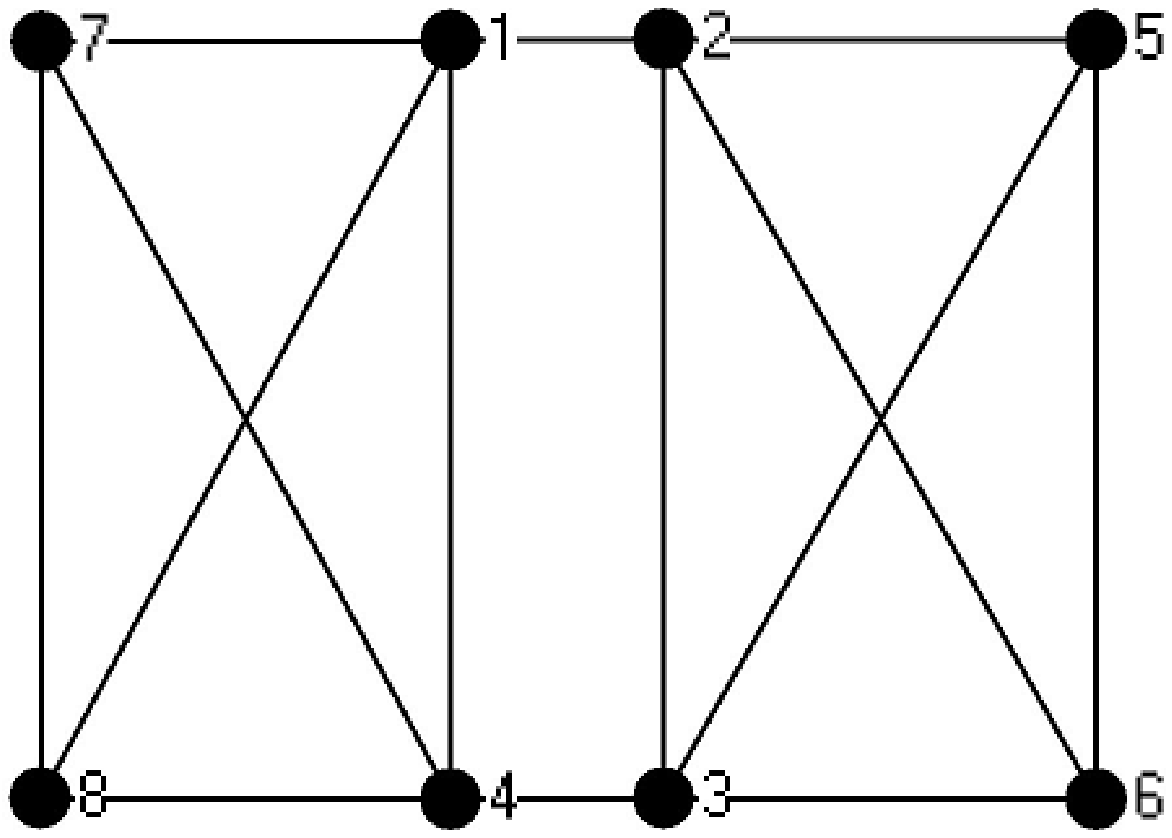}
\includegraphics[angle=0, width=3.5cm]{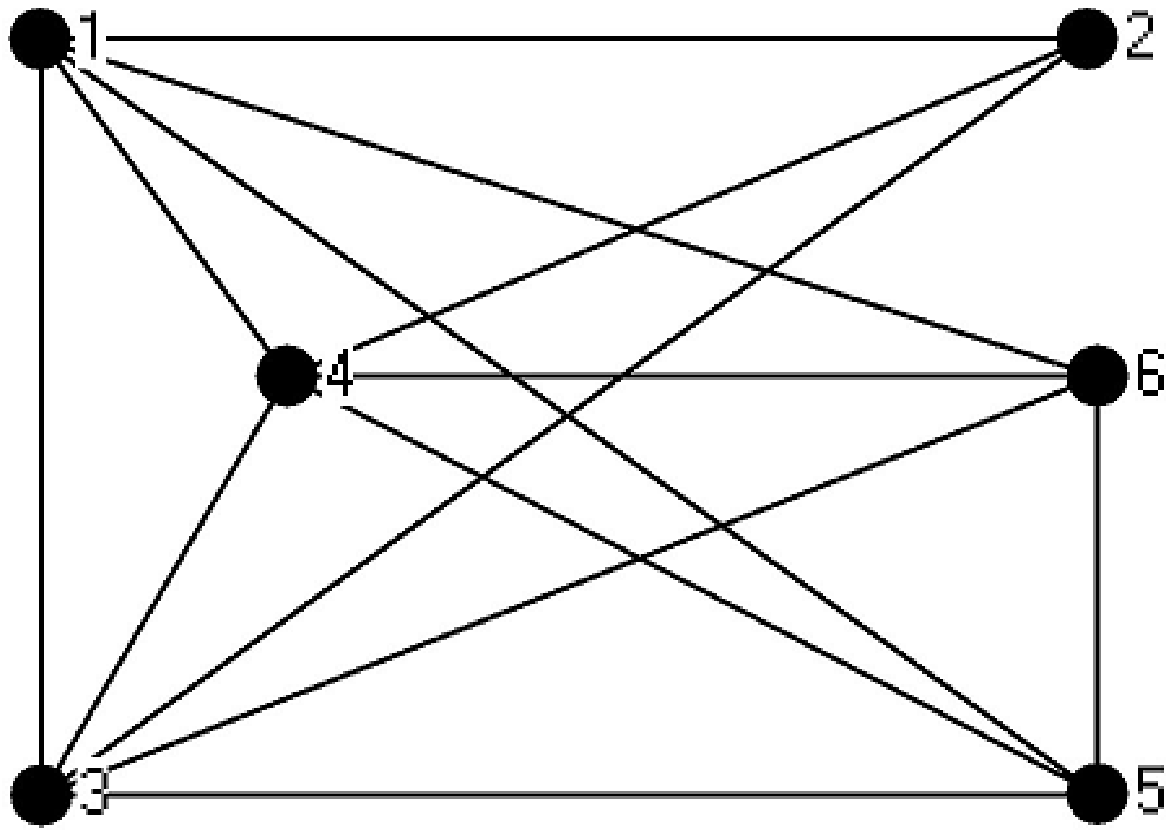}
\end{center}
%\vspace*{-1.8cm}
\end{figure}

The bounds {\it i}) and {\it ii}) are attained, for instance, for
the cocktail-party graph $\Gamma=K_6-F$ where $\gamma_2(\Gamma)=2$
and $\gamma_{\hat o}(\Gamma)=4$. The bound {\it iii}), is attained,
for instance, for the left hand side graph of Figure \ref{fig1}: in
this case $\gamma_{\hat o}(\Gamma)=6$. Example of equality in {\it
iv}) is $\Gamma=K_3\times K_2$. We emphasize that there are graphs
with minimum degree one, such that
 bounds {\it ii}) and {\it iii} fail. This is, for instance, the case of
the star graph, $\Gamma=K_{1,r}$, with $r\ge 6$. In this case
$n=r+1$ and $\gamma_{\hat{o}}(\Gamma)=\alpha(\Gamma)=r$.

\section{Bounding below the global offensive alliance number}

\begin{theorem}
For all connected graph $\Gamma$ of order $n$, minimum degree
$\delta$ and maximum degree $\Delta$,
\begin{itemize}
\item[i{\rm )}]  $
\gamma_0(\Gamma)\ge \left\lbrace \begin{array}{ll}
\left\lceil\frac{n(\delta+1)}{2\Delta+\delta+1}\right\rceil  & {\rm
if }\quad \delta \quad {\rm odd};
                            \\
                            \\
                            \left\lceil\frac{n\delta}{2\Delta+\delta}\right\rceil  &  {\rm otherwise;} \end{array}
                                                \right  .$

                                                \item[ii{\rm )}]  $
\gamma_{\hat{0}}(\Gamma)\ge \left\lbrace \begin{array}{ll}
\left\lceil\frac{n(\delta+3)}{2\Delta+\delta+3}\right\rceil  & {\rm
if }\quad \delta \quad {\rm odd};
                            \\
                            \\
                            \left\lceil\frac{n(\delta+2)}{2\Delta+\delta+2}\right\rceil  &  {\rm otherwise.} \end{array}
                                                \right  .$
\end{itemize}
\end{theorem}

\begin{proof}
Let $\gamma_{_k}(\Gamma)$ denotes the $k$-domination number of
$\Gamma$. Since all global strong offensive alliance is a
$\left\lceil\frac{\delta+1}{2}\right\rceil$-dominating set and all
global strong offensive alliance is a
$\left\lceil\frac{\delta+2}{2}\right\rceil$-dominating set,
\begin{equation}
\gamma_{_{\left\lceil\frac{\delta+1}{2}\right\rceil}}(\Gamma) \leq
\gamma_{0}(\Gamma)
\end{equation} and
\begin{equation}
\gamma_{_{\left\lceil\frac{\delta+2}{2}\right\rceil}}(\Gamma) \leq
\gamma_{\hat{0}}(\Gamma).
\end{equation}
 On the other hand, for all k-dominating
set $S\subset V$, $k(n-|S|)\leq \Delta |S|.$ Hence,
\begin{equation} \label{k-domination}
\gamma_{_k}(\Gamma)\ge \left\lceil\frac{kn}{\Delta+k}\right\rceil.
\end{equation}
Therefore, the result follows.
\end{proof}

 Examples of equality in above theorem are
 $\Gamma=K_{3,3}$ and the $3$-cube graph.

The following result provides  tight bounds on $\gamma_{o}(\Gamma)$
and $\gamma_{\hat{o}}(\Gamma)$ in terms of the order and size  of
$\Gamma$.

\begin{theorem}\label{Thcotainf}
For all graph  $\Gamma$  of order  $n$ and size $m$,
$$\gamma_{o}(\Gamma)\ge \left\lceil\frac{3n-\sqrt{9n^{2}-8n-16m}}{4}\right\rceil$$ and
$$\gamma_{\hat{o}}(\Gamma)\ge \left\lceil\frac{3n+1-\sqrt{9n^2-10n-16m+1}}{4}\right\rceil.$$
\end{theorem}

\begin{proof}
If $S$ denotes a global offensive alliance in  $\Gamma=(V,E)$, then

\begin{equation}\label{Global11}
2m=\sum_{v\in V\backslash S}\delta(v)+\sum_{v\in S}\delta(v) \le
(n-|S|)(2|S|-1)+|S|(n-1).
\end{equation}
Hence, solving $2|S|^{2}-3n|S|+2m+n\leq 0$ we obtain the bound on
$\gamma_{o}(\Gamma)$. The bound on $\gamma_{\hat{o}}(\Gamma)$ is
derived by analogy by using
\begin{equation}\label{Global111}
2m\leq (n-|S|)(2|S|-2)+|S|(n-1) \end{equation} instead of
(\ref{Global11}).
\end{proof}

Example of equality in the above bounds is the right hand side graph
of Figure \ref{fig1} where $S=\{2,6,5\}$ is a minimal global
offensive alliance and $S'=\{1,3,4\}$ is a minimal  global strong
offensive alliance. Even so, the following bounds, expressed  in
terms of the order, size, and the maximum degree of $\Gamma$,
improve the previous result.

\begin{theorem}
For all graph  $\Gamma$  of order  $n$, size $m$ and maximum degree
$\Delta$,
$$\gamma_{0}(\Gamma)\ge
\left\lceil\frac{2m+n}{3\Delta+1}\right\rceil \quad {\rm  and} \quad
\gamma_{\hat{0}}(\Gamma)\ge
\left\lceil\frac{2(m+n)}{3\Delta+2}\right\rceil.$$
\end{theorem}

\begin{proof}
 If $S\subset V$, then
\begin{equation}\label{ll1}
|S|\Delta \ge \sum_{v\in V\backslash S}|N_{ S}(v)|.
\end{equation}
Moreover, if $S$ is a global offensive alliance in  $\Gamma$, then
\begin{equation}\label{Globall1}
\sum_{v\in V\backslash S}|N_{S}(v)|\geq\sum_{v\in V\backslash
S}|N_{V\backslash S}(v)|+(n-|S|).
\end{equation}
Thus,
\begin{align*}
2m&=\sum_{v\in V\backslash S}\delta(v)+\sum_{v\in S}\delta(v) \\
  &=
\sum_{v\in V\backslash S}|N_S(v)|+\sum_{v\in V\backslash
S}|N_{V\backslash
S}(v)|+\sum_{v\in S}\delta(v) \\
& \le 2\sum_{v\in V\backslash S}|N_S(v)|+|S|-n+\sum_{v\in S}\delta(v) \\
& \le (3\Delta +1)|S|-n.
\end{align*}
So, the  bound on $\gamma_{0}(\Gamma)$ follows. If the global
offensive alliance $S$ is strong, then we have
\begin{equation}\label{f}
\sum_{v\in V\backslash S}|N_{S}(v)|\ge \sum_{v\in V\backslash
S}|N_{V\backslash S}(v)|+2(n-|S|).
\end{equation}
Basically, the bound on $\gamma_{\hat{o}}(\Gamma)$ follows as before
by using (\ref{f}) instead of (\ref{Globall1}).
\end{proof}

The above bounds are reached, for instance, in the case of the
3-cube graph $\Gamma=K_2\times K_2\times K_2$, where
$\gamma_o(\Gamma)=\gamma_{\hat{o}}(\Gamma)=4$. Notice that Theorem
\ref{Thcotainf} only gives $\gamma_o(\Gamma)\ge 2$.

As we can see in \cite{spectral}, we can obtain bounds on the
alliance numbers from the spectrum of $\Gamma$ or from the Laplacian
spectrum of $\Gamma$. For instance, the following result was proved
in \cite{spectral}. For completeness we include the proof of this
result.

\begin{theorem}\label{thOf}
For all graph  $\Gamma$  of order  $n$ and size $m$, minimum degree
$\delta$ and  Laplacian spectral radius $\mu$,
$$\gamma_{o}(\Gamma)\ge \left\lceil\frac{n}{\mu}\left\lceil\frac{\delta +1}{2}\right\rceil\right\rceil
\quad {\rm and } \quad  \gamma_{\hat{o}}(\Gamma)\ge
\left\lceil\frac{n}{\mu}\left(\left\lceil\frac{\delta
}{2}\right\rceil+1\right)\right\rceil.$$
\end{theorem}

\begin{proof}

It was shown in  \cite{fiedler} that the Laplacian spectral radius
of $\Gamma$, $\mu$, satisfies
\begin{equation}\label{rfiedler1}
  \mu=2n \max \left\{ \frac{\sum_{v_i\sim v_j}(w_i-w_j)^2  }
  {\sum_{v_i\in V}\sum_{v_j\in V}(w_i-w_j)^2}: \mbox{\rm $w\neq \alpha{\bf j}$
  for   $\alpha\in \mathbb{R}$ } \right\}.
\end{equation}
Let $S\subset V$.  From (\ref{rfiedler1}), taking $w\in
\mathbb{R}^n$ defined as
$$
w_i= \left\lbrace \begin{array}{ll} 1  & {\rm if }\quad  v_i\in S;
                            \\ 0 &  {\rm otherwise,} \end{array}
                                                \right  .$$
                                                we obtain
\begin{equation}\label{fiedlerAlliance1}
\mu\ge \frac{n\displaystyle\sum_{v\in V\setminus S} | N_{ S}(v) |
}{|S|(n-|S|)}.
\end{equation}
Moreover, if $S$ is a global offensive alliance  in $\Gamma$,
\begin{equation} \label{globalGrado1}
|N_{ S}(v)|\ge \left\lceil\frac{\delta(v)+1}{2}\right\rceil \quad
\forall v\in V\setminus S.
\end{equation}
 Thus,
  (\ref{fiedlerAlliance1}) and (\ref{globalGrado1}) lead to
\begin{equation}\label{final1}
\mu\ge \frac{n}{|s|}\left\lceil\frac{ \delta  +1}{2}\right\rceil.
\end{equation}
Therefore, solving (\ref{final1}) for $|S|$, and considering that it
is an integer, we obtain the bound on $\gamma_{{{a}_o}}(\Gamma)$. If
the  global offensive alliance $S$ is strong, then
\begin{equation} \label{globalGrado2}
|N_{ S}(v)|\ge \left\lceil\frac{\delta(v)}{2}\right\rceil +1\quad
\forall v\in V\setminus S.
\end{equation}
Thus,
  (\ref{fiedlerAlliance1}) and (\ref{globalGrado2}) lead to  the bound on
  $\gamma_{\hat{o}}(\Gamma)$.
\end{proof}

If $\Gamma$ is the Petersen graph, then $\mu=5$. Thus, Theorem
\ref{thOf} leads to $\gamma_{o}(\Gamma)\ge 4$ and
$\gamma_{\hat{o}}(\Gamma)\ge 6.$ Therefore, the above bounds are
tight.

\section{Offensive alliances and connected subgraphs}

An offensive alliance (global offensive alliance) $S$ in $\Gamma$ is
\emph{minimal} if no proper subset of $S$ is an offensive alliance
(global offensive alliance) in $\Gamma$.

\begin{theorem} \label{diametofen}
Let  $\Gamma=(V,E)$  be a  connected graph of order $n$ and diameter
$D(\Gamma)$. If $\Gamma$ has a minimal (global) offensive alliance
$S$ such that $\langle V\backslash S \rangle$ is connected, then
$D(\Gamma)\leq n-|S|+1.$
\end{theorem}

\begin{proof}
If  $S\subset V$ is a minimal (global) offensive alliance  in
$\Gamma$ then  $ V\backslash S$ is a dominating set in $\Gamma$. So,
if $\langle V\backslash S \rangle$ is connected, then $D(\Gamma)\leq
D(\langle V\backslash S \rangle)+2.$ Hence, $D(\Gamma)\leq n-|S|+1.$
\end{proof}

We remark that there are graphs such that for every minimal (global)
offensive alliance $S$, $\langle V \backslash S\rangle$ is not
connected. For instance,  the case of the 3-cube graph.

\begin{figure}[h]
\begin{center}
\caption{} \label{fig2} %\vspace{-1,5cm}
\includegraphics[angle=0, width=4.5cm]{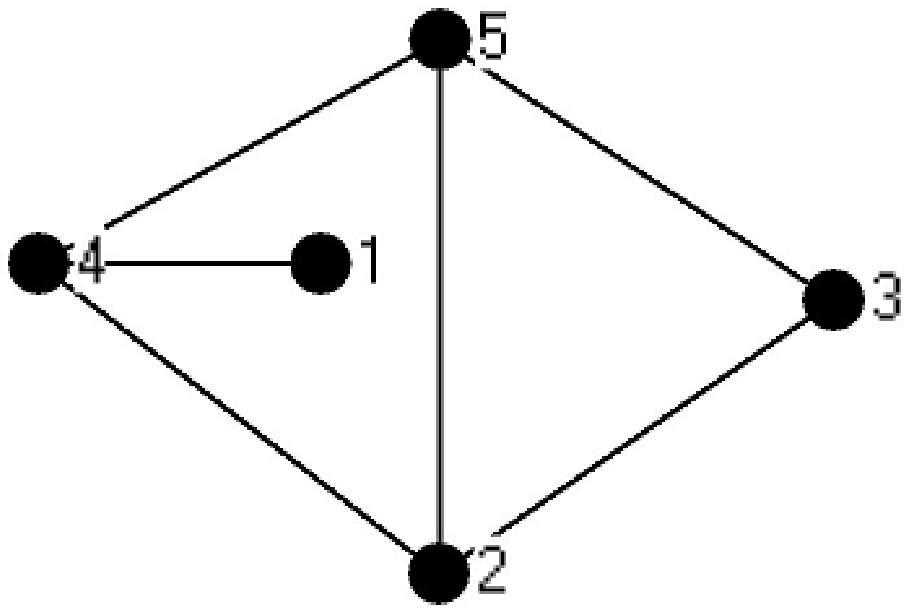}
\includegraphics[angle=0, width=4.5cm]{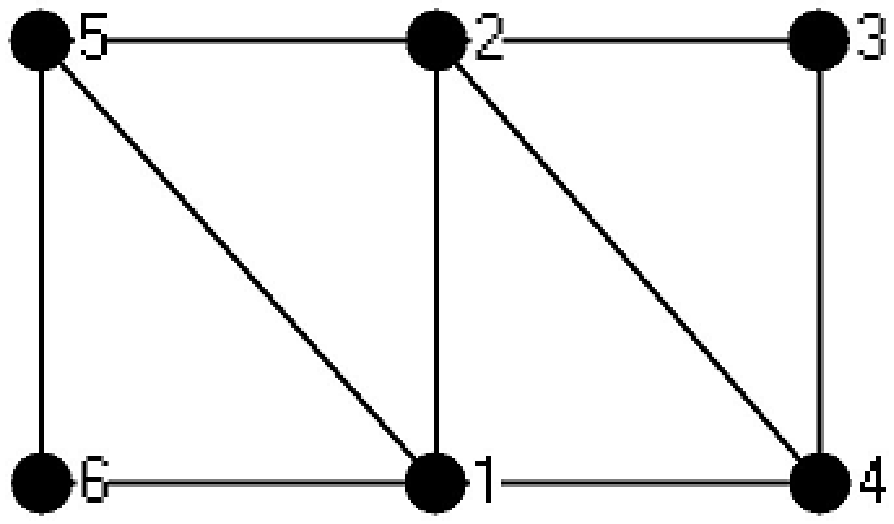}
\end{center}
%\vspace*{-1.8cm}
\end{figure}

The above bound is tight.  Let $\Gamma$ be the left hand side graph
of Figure  \ref{fig2}. In this case the set  $S=\{1,3,5\}$ is a
minimal global offensive alliance and $ V\backslash S=\{2,4\}$ is
connected. Thus,  $3=D(\Gamma)\leq n-|S|+1=3$.

\begin{theorem}
Let  $\Gamma=(V,E)$  be a graph of order  $n$ and maximum degree
$\Delta$. For all minimal global offensive alliance $S$ such that
$\langle V\backslash S \rangle$ is connected,
$$|S|\ge \left\lceil\frac{3n-2}{\Delta+3}\right\rceil.$$
Moreover, for all minimal global strong offensive alliance $S$ such
that $\langle V\backslash S \rangle$ is connected,
$$|S|\ge \left \lceil\frac{4n-2}{\Delta+4}\right\rceil.$$
\end{theorem}

\begin{proof}
Let $S\subset V$. As $\langle V\backslash S \rangle$ is connected,
\begin{equation}\label{l2}
 \sum_{v\in V\backslash
S}|N_{V\backslash S}(v)|\geq 2(n-|S|-1).
\end{equation}
So, the first  bound follows, by (\ref{ll1}), (\ref{Globall1}) and
(\ref{l2}).  The second bound is derived by analogy by using
(\ref{f}) instead of (\ref{Globall1}).
\end{proof}

%XXXXXXXXXXXXXXXXX ofensiva XXXXXXXXXXXXXXXXXXXx
%\begin{figure}[h]
%\begin{center}
%\caption{ } \label{fig2}
%\includegraphics[width=0.3\textwidth]{ofensivadiametro}
%\includegraphics[width=0.3\textwidth]{ofenteorema8b}
%\end{center}
%\end{figure}

The above bounds are tight. If $\Gamma$ is the left hand side graph
of Figure \ref{fig2}, then $S=\{1,3,5\}$ is a minimal global
offensive alliance in $\Gamma$ and $V\backslash S=\{2,4\}$ is
connected. Moreover, if $\Gamma$ is the right hand side graph of
Figure \ref{fig2}, then $S=\{3,4,5,6\}$ is a minimal global strong
offensive alliance in $\Gamma$ and $V\backslash S=\{1,2\}$ is
connected.

We define the {\em global-connected offensive allian\-ce number},
 $\gamma_{co}(\Gamma)$,
 (respectively, {\em global-connected strong offensive
alliance number} $\gamma_{c\hat{o}}(\Gamma)$) as the minimum
cardinality of any global offensive alliance (respectively, global
strong offensive alliance) in $\Gamma$ whose induced subgraph is
connected.

\begin{theorem}\label{ThOffConnected}
Let $\Gamma$ be a simple graph of order $n$, size $m$, diameter $D$
and maximum degree $\Delta$. The global-connected offensive alliance
number of $\Gamma$ is bounded by
$$\gamma_{co}(\Gamma)\ge \left\lceil\frac{2m+n+2(D-1)^2}{2n+\Delta+1}\right\rceil$$
and the global-connected strong offensive alliance number of
$\Gamma$ is bounded by
$$\gamma_{\hat{co}}(\Gamma)\ge \left\lceil\frac{2\left(m+n+(D-1)^2\right)}{2n+\Delta+2}\right\rceil.$$
\end{theorem}

\begin{proof}
If  $S$ is a  global offensive alliance in $\Gamma=(V,E)$, then by
(\ref{Globall1}) we have
\begin{equation}\label{ofensivaGlobal1}
(|S|-1)(n-|S|)\ge \sum_{v\in V\backslash S}|N_{V\backslash S}(v)|.
\end{equation}
Thus,
\begin{equation}\label{d}
 (2|S|-1)(n-|S|)\ge \sum_{v\in V\backslash S}|N_{S}(v)|+
\sum_{v\in V\backslash S}|N_{V\backslash S}(v)|=\sum_{v\in
V\backslash S}\delta(v).
\end{equation}
Therefore,
\begin{equation}\label{e}
(2|S|-1)(n-|S|)+\Delta|S|\ge\sum_{v\in V\backslash
S}\delta(v)+\sum_{v\in S}\delta(v)=2m.
\end{equation}
On the other hand, if $S$ is a dominating set and  $\langle S
\rangle$ is connected, then $D(\Gamma)\leq D(\langle S \rangle)+2.$
So, $D(\Gamma)\le |S|+1$. Hence,
\begin{equation}\label{condiam}
2n|S|-n+|S|+\Delta|S|\ge 2m+2(D(\Gamma)-1)^2.
\end{equation}
Thus, the bound on $\gamma_{co}(\Gamma)$ follows. Basically the
bound on $\gamma_{\hat{co}}(\Gamma)$ follows as before by using
(\ref{f}) instead of (\ref{Globall1}).
\end{proof}

The above bounds are tight, as we show in the following instance.
Let $\Gamma_{3,t}$ be the graph obtained by joining every vertex of
the complete graph $K_3$ with every vertex of the trivial graph of
order $t\ge 8$. In such case,
$\gamma_{co}(\Gamma_{3,t})=\gamma_{\hat{co}}(\Gamma_{3,t})=3$ and
Theorem \ref{ThOffConnected} leads to $\gamma_{co}(\Gamma_{r,t})\ge
3$ and $\gamma_{\hat{co}}(\Gamma_{3,t})\ge 3$.


\begin{thebibliography}{99}


\bibitem{Cockayne} E. J. Cockayne, B. Gamble,  B. Shepherd,  An upper bound for the
$k$-domination number of a graph. J. Graph Theory {\bf 9} (4) (1985)
533-534.



\bibitem{favaron}  O. Favaron, G.  Fricke,  W. Goddard, S. M. Hedetniemi,  S. T. Hedetniemi,
 P. Kristiansen, R. C. Laskar and D. R. Skaggs,
Offensive alliances in graphs. {\it Discuss. Math. Graph Theory }
{\bf 24} (2)(2004), 263-275.



\bibitem{fiedler}  M. Fiedler,  A property of eigenvectors of nonnegative symmetric
  matrices and its application to graph theory,
  {\it  Czechoslovak Math. J.} {\bf  25} (100) (1975), 619-633.

\bibitem{GlobalalliancesOne} T. W. Haynes,  S. T. Hedetniemi,  and M. A. Henning,
Global defensive alliances in graphs, {\it Electron. J. Combin.}
{\bf 10} (2003), Research Paper 47.



\bibitem{cconex} S.T. Hedetniemi, R. Laskar, Connected domination in graphs, Graph Theory and Combinatorics:
Proceedings of the Cambridge Combinatorial Conference, Academic
Press, London, 1984.


\bibitem{alliancesOne} P. Kristiansen, S. M. Hedetniemi and S. T.
Hedetniemi, Alliances in graphs. {\it J. Combin. Math. Combin.
Comput.} {\bf 48} (2004), 157-177.

%\bibitem{li} L. Jiong-Sheng and Z. Xiao-Dong, A New Upper bound for eigenvalues
%of the Laplacian Marrix of a Graph. {\it Linear Algebra and
%Applications.} {\bf 265}(1997) 93-100.


\bibitem{spectral} J. A. Rodr\'{\i}guez and J. M. Sigarreta,
Spectral study of alliances in graphs.  Submitted 2005.

\bibitem{partition} J. A. Rodr\'{\i}guez, Laplacian eigenvalues and partition problems in hypergraphs.
\emph{Math. Preprint Archive}, {\bf  2004}, Issue 3 (2004) 183-196.
 \emph{Linear Algebra and its Applications.}  Submitted 2003.



\bibitem{planar} J. A. Rodr\'{\i}guez and J. M. Sigarreta,  Global alliances in planar
graphs. Submitted 2005.

\bibitem{linedefensive}  J. M. Sigarreta and J. A. Rodr\'{\i}guez,  On  defensive alliances and line
graphs. \emph{Applied Mathematics Letters}. In press.



\end{thebibliography}
\end{document}